\documentclass{gtart_h}

\def\ifplaintex{\expandafter\ifx\csname documentclass\endcsname\relax}

\def\gtp{{\mathsurround=0pt\it $\cal G\mskip-2mu$eometry \&\ 
$\cal T\!\!$opology $\cal P\!$ublications}}  

\def\recd{{\small Received:\qua\receiveddate\ifx\reviseddate\relax
\else\qquad Revised:\qua\reviseddate\fi\par}} 

\def\lognumber#1{\def\thelognumber{#1}}
\def\volumenumber#1{\def\thevolumenumber{#1}}
\def\volumeyear#1{\def\thevolumeyear{#1}}
\def\papernumber#1{\def\thepapernumber{#1}}
\def\pagenumbers#1#2{\def\startpage{#1}\def\finishpage{#2}}
\def\published#1{\def\publishdate{#1}}

\def\received#1{\def\receiveddate{#1}}

\def\accepted#1{\def\accepteddate{#1}}

\def\asciiaddress#1{\def\theasciiaddress{#1}}
\def\asciiemail#1{\def\theasciiemail{#1}}

\long\def\asciiabstract#1{\long\def\theasciiabstract{#1}}

\let\\\par\let\thelognumber\relax\let\thevolumenumber\relax
\let\thepapernumber\relax\let\thevolumeyear\relax\let\startpage\relax
\let\finishpage\relax\let\publishdate\relax\let\receiveddate\relax
\let\reviseddate\relax\let\accepteddate\relax\let\theasciititle\relax
\let\theasciiauthors\relax\let\theasciiaddress\relax
\let\theasciiabstract\relax

\let\theasciiemail\relax

\ifplaintex
\font\logobig=cmssbx10 scaled 3836
\font\logomed=cmssbx10 scaled 2557
\else
\font\logobig=cmssbx10 scaled 4200
\font\logomed=cmssbx10 scaled 2800
\fi

\long\def\makeagttitle{   
\count0=\startpage
\agt\hfill      
\hbox to 45truept{\vbox to 0pt{\vglue -13truept{\logomed A\kern -.37em{\logobig 
T}\kern -.38em G}\vss}\hss}
\break
{\small Volume \thevolumenumber\ (\thevolumeyear)
\startpage--\finishpage\nl
Published: \publishdate}

\vglue .25truein

{\parskip=0pt\leftskip 0pt plus
1fil\def\\{\par\smallskip}{\Large\bf\thetitle}\par\medskip} \vglue
0.05truein

{\parskip=0pt\leftskip 0pt plus 1fil\def\\{\par}{\sc\theauthors}
\par\medskip}%
 
\vglue 0.03truein 

{\small\leftskip 25truept\rightskip 25truept{\bf Abstract}\stdspace\theabstract

{\bf AMS Classification}\stdspace\theprimaryclass
\ifx\thesecondaryclass\relax\else; \thesecondaryclass\fi\par
{\bf Keywords}\stdspace \thekeywords\par}\vglue 7truept

}   

\ifplaintex
\font\phead=cmsl9 scaled 950
\font\pnum=cmbx10 scaled 913
\font\pfoot=cmsl9 scaled 950
\headline{\vbox to 0pt{\vskip -4.5mm\line{\small\phead\ifnum
\count0=\startpage ISSN 1472-2739 (on-line) 1472-2747 (printed)
\hfill {\pnum\folio}\else\ifodd\count0\def\\{ }%
\ifx\theshorttitle\relax\thetitle\else\theshorttitle\fi\hfill{\pnum\folio}
\else\def\\{ and }{\pnum\folio}\hfill\ifx\theshortauthors\relax\theauthors
\else\theshortauthors\fi\fi\fi}\vss}}
\footline{\vbox to 0pt{\vglue 0mm\line{\small\pfoot\ifnum\count0=\startpage
\copyright\ \gtp\hfill\else
\agt, Volume \thevolumenumber\ (\thevolumeyear)\hfill\fi}\vss}}
\else
\font=cmsl9 scaled 1050
\font\lnum=cmbx10 
\font=cmsl9 scaled 1050
\makeatletter
\def\@oddhead{{\small\lhead\ifnum\count0=\startpage ISSN 1472-2739 
(on-line) 1472-2747 (printed)\hfill {\lnum\number\count0}\else\ifodd\count0
\def\\{ }\ifx\theshorttitle\relax \thetitle \else\theshorttitle\fi\hfill
{\lnum\number\count0}\else\def\\{ and }{\lnum\number\count0}
\hfill\ifx\theshortauthors\relax 
\theauthors\else\theshortauthors\fi\fi\fi}}\def\@evenhead{\@oddhead}
\def\@oddfoot{\small\lfoot\ifnum\count0=\startpage\copyright\ \gtp\hfill\else
\agt, Volume \thevolumenumber\ (\thevolumeyear)\hfill\fi}
\def\@evenfoot{\@oddfoot}
\makeatother
\fi
\let\maketitlepage\makeagttitle
\let\makeshorttitle\maketitlepage
\let\maketitle\maketitlepage

\newwrite\gtoutfile
\long\gdef\makeheadfile{  
{\def\\{, }\def\s{ }
\immediate\openout\gtoutfile head.xxx
\immediate\write\gtoutfile{Proxy-for: \ifx\theasciiauthors\relax
\theauthors\else\theasciiauthors\fi\s<\ifx\theasciiemail\relax\theemail\else\theasciiemail\fi>}
\immediate\write\gtoutfile{\noexpand\\}
\immediate\write\gtoutfile{Authors: \ifx\theasciiauthors\relax
\theauthors\else\theasciiauthors\fi}
{\def\\{ }\immediate\write\gtoutfile{Title: \ifx\theasciititle\relax
\thetitle\else\theasciititle\fi}}
\immediate\write\gtoutfile{Subj-class: GT or SG, GR etc}
\immediate\write\gtoutfile{MSC-class: \theprimaryclass\ifx\thesecondaryclass\relax\else, \thesecondaryclass\fi}
\immediate\write\gtoutfile{Journal-ref: Algebr. Geom. Topol. \thevolumenumber\s
(\thevolumeyear) \startpage-\finishpage}
\immediate\write\gtoutfile{Comments: Published by Algebraic and
Geometric Topology at}
\immediate\write\gtoutfile{\s\s\s  http://www.maths.warwick.ac.uk/agt/AGTVol\thevolumenumber/agt-\thevolumenumber-\thepapernumber.abs.html}
\immediate\write\gtoutfile{\noexpand\\}
\immediate\write\gtoutfile{}
\ifx\theasciiabstract\relax
\immediate\write\gtoutfile{\theabstract}\else
\immediate\write\gtoutfile{\theasciiabstract}\fi
\immediate\write\gtoutfile{}
\immediate\write\gtoutfile{\noexpand\\}
\immediate\write\gtoutfile{}
\immediate\closeout\gtoutfile}}  

\def\maketitlepage{\makeagttitle\makeheadfile}
\let\makeshorttitle\maketitlepage
\let\maketitle\maketitlepage

\lognumber{22}
\volumenumber{5}
\volumeyear{2005}
\papernumber{22}
\pagenumbers{509}{535}
\received{31 May 2004} 
\accepted{16 April 2005}
\published{4 June 2005}

\usepackage{amsmath,amssymb,epsf}
\allowdisplaybreaks

\newtheorem{thm}{Theorem}[section]
\newtheorem{cor}[thm]{Corollary}
\newtheorem{lem}[thm]{Lemma}
\newtheorem{prop}[thm]{Proposition}

\theoremstyle{remark}

\newtheorem*{rem}{Remark}

\begin{document}

\title{Minimal surface representations\\of virtual knots and links}

\author{H.A. Dye\\Louis H. Kauffman}

\address{MADN-MATH, United States Military Academy\\
646 Swift Road, West Point, NY 10996, USA\\{\rm 
and}\\Department of Mathematics, Statistics and Computer 
Science\\University of Illinois at Chicago, 851 South Morgan St\\Chicago, 
IL  60607-7045, USA}

\asciiaddress{MADN-MATH, United States Military Academy\\
646 Swift Road, West Point, NY 10996, USA\\
and\\Department of Mathematics, Statistics and Computer 
Science\\University of Illinois at Chicago, 851 South Morgan St\\Chicago, 
IL  60607-7045, USA}

\gtemail{\mailto{hdye@ttocs.org}, \mailto{kauffman@uic.edu}}
\asciiemail{hdye@ttocs.org, kauffman@uic.edu }


\begin{abstract}
Kuperberg \cite{kup} has shown that a virtual knot diagram corresponds
(up to generalized Reidemeister moves) to a unique embedding in a
thickened surface of minimal genus.  If a virtual knot diagram is
equivalent to a classical knot diagram then this minimal surface is a
sphere.  Using this result and a generalised bracket polynomial, we
develop methods that may determine whether a virtual knot diagram is
non-classical (and hence non-trivial).  As examples we show that,
except for special cases, link diagrams with a single virtualization
and link diagrams with a single virtual crossing are non-classical.
\end{abstract}

\asciiabstract{%
Kuperberg [Algebr. Geom. Topol. 3 (2003) 587-591] has shown that a
virtual knot corresponds (up to generalized Reidemeister moves) to a
unique embedding in a thickened surface of minimal genus.  If a
virtual knot diagram is equivalent to a classical knot diagram then
this minimal surface is a sphere.  Using this result and a generalised
bracket polynomial, we develop methods that may determine whether a
virtual knot diagram is non-classical (and hence non-trivial).  As
examples we show that, except for special cases, link diagrams with a
single virtualization and link diagrams with a single virtual crossing
are non-classical.}

\maketitle

\section{Introduction}
Virtual knot diagrams are a generalization of classical knot diagrams
introduced by L. Kauffman in 1996 \cite{kvirt}. Results in this area
immediately indicated that the bracket polynomial and the fundamental
group did not detect many non-trivial and non-classical virtual knot
diagrams. We are interested in detecting non-trivial virtual knot
diagrams, and, in particular, determining if a virtual knot diagram is
non-classical and non-trivial.

The bracket polynomial and the fundamental group can not differentiate
all non-trivial virtual knot diagrams from the unknot. Kauffman, in
\cite{kvirt}, gave a process of virtualization that produces from a
diagram $ K$, a pair of diagrams: a virtual knot diagram $K_v $ and a
classical knot diagram $ K_s $ (obtained by switching a crossing in
$K$). The diagrams $K_v $ and $ K_s $ have the same bracket
polynomial. One can show that if $ K $ is a non-trivial classical
knot, then $ K_v $ is a non-trivial virtual knot (possibly classical)
\cite{kvirt}.  This process may be used to construct non-trivial
virtual knot diagrams with trivial bracket polynomial. There are also
other virtual knot diagrams with trivial bracket polynomial that are
not produced by virtualization.  Kishino's knot is the first example
of this type and it is not differentiated from the unknot by the
fundamental group or the Jones polynomial.  In \cite{kishpoly},
Kishino's knot was detected by the 3-strand bracket
polynomial. Kishino's knot is also detected by the quaternionic
biquandle \cite{BIQ}, \cite{kmant}.  However, these invariants can be
difficult to compute.  Other virtual knot diagrams, undetected by the
fundamental group and the bracket polynomial, are described in
\cite{dye}. It is shown in \cite{dye} that there are an infinite
number of virtual knot diagrams that are not detected by the
fundamental group or the Jones polynomial.

Using the bracket polynomial and Kuperberg's result \cite{kup}, we
develop methods that may determine if a virtual knot diagram is
non-trivial and non-classical. We focus on the case of knot diagrams
with one virtualization and the examples in \cite{dye}.  We show that,
except for special cases, link diagrams with a single virtualization
and link diagrams with a single virtual crossing are non-classical and
non-trivial. We construct examples of virtual link diagrams with
either one virtualization or one virtual crossing using the methods
from \cite{morwen} that are not detectable by the surface bracket
polynomial.  In the final section, we discuss virtual knots produced
by two virtualizations.

\section{Virtual knots and links and minimal representations}

Virtual knot (and link) theory was introduced by the second author in
\cite{kvirt}, to which we refer for basic concepts and notation.  In
particular we use a small circle to indicate a virtual crossing in a
diagram.  A virtual knot (or link) is an equivalence class of diagrams
containing ordinary crossings and virtual crossings under the three
familar Reidemeister moves together with the ability to move an arc
containing only virtual crossings to any other position with the same
endpoints.  A virtual link diagram can be \emph{represented\/} as a
link diagram in an oriented surface, by adding a handles to
desingularize the virtual crossings.  This can also be done by
regarding the diagram as lying in a surface with boundary (cf Kamada
and Kamada \cite{kamada}) and capping the boundary components.  This
surface can be \emph{stabilized\/} by adding further handles which do
not meet the diagram.  This diagram can then be regarded as a genuine
link in the surface thickened by crossing with a unit interval.
Carter, Kamada and Saito \cite{cks} have shown that desingularization
induces a bijection between equivalence classes of virtual links and
stable equivalence classes of links in a thickened surface.  A short
proof of this result and a good summary of the various different ways
of thinking of a virtual link can be found in \cite[Theorem 4.5 and
above]{racks}, see also \cite{Kdetect,Keldysh}.  In \cite{Keldysh} the
Kuperberg result (below) is used to show that virtual knots are
algorithmically recognizable and to find the genus of connected sums
of virtual knots.

A representation of a surface is \emph{minimal\/} if it cannot be
destabilized.   

\begin{thm}[Kuperberg \cite{kup}]A minimal representation of a virtual
link is unique up to homeomorphism of the thickened surface.
\end{thm}

\begin{cor}\label{mingenus}If the minimal surface has genus greater 
than zero then the link is non-trivial and not equivalent to a
classical link.
\end{cor}

We need to think of these results in terms of diagrams, rather than
thickened surfaces.  In these terms destabilization and is performed
by surgering the surface along a \emph{cancellation curve} which does
not meet the diagram.  A representation of a link is then minimal if
no cancellation curve can be found after a possible sequence of
Reidemeister moves.  We refer to a minimal representation as a
\emph{characterization} and we define the \emph{virtual genus} of a
virtual knot or link diagram to be the minimal genus.

Corollary \ref{mingenus} taken alone does not provide an algorithm
that determines cancellation curves.  However, such algorithms can be
formulated using normal surface theory \cite{kup2}, \cite{nst}.  In
this paper we use algebraic invariants of knots and links to
investigate minimality.

We need to stress the non-constructive nature of minimality:

\begin{rem}\label{virtgen}To obtain a characterization from a
representation it may be necessary to perform a sequence of handle
cancellations and Reidemeister moves combined with non-trival
homeomorphisms of the surfaces.  \end{rem}

We will introduce new methods which apply a generalization of the
bracket polynomial to representations of virtual knot diagrams.  These
methods often determine if a given representation is minimal and
hence, by Corollary \ref{mingenus}, if a virtual knot diagram is
non-classical and non-trivial.  One method uses homology classes and
intersection numbers to determine non-triviality and the other method
uses isotopy classes to determine non-triviality.  Both methods
utilize the bracket skein relation to produce states that consist of
simple closed curves in the surface with coefficients in $ \mathbb{Z}
[A,A^{-1}] $ from a fixed representation. We make the following
definitions.

For a fixed representation of a virtual knot diagram, we refer to 
the \emph{surface-knot pair}, $(F,K)$, to indicate a specific choice of
surface and embedding of the knot.  
A \emph{surface-state pair}, $(F,s)$,
is a collection of disjoint simple closed curves in the surface.

We obtain a surface-state pair $(F,s)$ from $(F,K)$ by assigning a smoothing 
type to each classical crossing in the surface. 

\begin{figure}[htb] \epsfysize = 0.8 in
\centerline{\epsffile{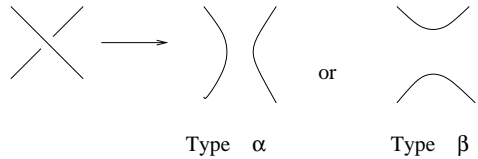}}
\caption{Smoothing Types}
\label{fig:smooth}
\end{figure}

Note that if $(F,K)$ has $ n $ classical crossings we 
obtain $ 2^n $ surface-state pairs, denoted 
$ \lbrace (F,s_1), ... (F, s_{2^n}) \rbrace $, by assignment of smoothing type. We denote the collection of all surface-state pairs as $ (F,S) $.
A surface state pair is analogous to a state of a classical knot diagram.  

We consider a surface-state pair, $ (F, s) $, in more detail. Each surface-state 
pair is a collection of disjoint
simple closed curves on the surface $ F $.
For a fixed collection of disjoint sets of curves in the surface F
we may study either isotopy classes of these curves (with or without 
homeomorphisms) or homology classes of curves.

We define the \emph{surface bracket polynomial} of a 
representation  $ (F,K) $.
 Let $ \hat{K} $ be a
 virtual knot diagram 
and let  $ (F, K ) $ be a fixed representation of $ \hat{K} $. 
The surface bracket polynomial of $ K $ 
is denoted as $ \langle (F,K) \rangle $. Then:
$$ 
 \langle (F,K) \rangle = \underset{(F,s(c)) \in (F,S)}{ \sum} \langle K| s(c) \rangle d^{|s(c)|} [s(c)]
$$
where $ \langle K | s(c) \rangle = A^{c(s)} $ and $ c(s) $ is the number of type $ \alpha
$ smoothings minus the number of type $ \beta $ smoothings. 
$ | s(c) | $ is the number of curves which bound a disk in the surface 
and $ [ s(c) ] $ represents a formal sum
 of the disjoint curves that do not bound a disk in 
the surface-state pair $ (F,s(c)) $.

(Note that $ [s(c)] $ may be regarded either as a formal sum of
homology classes in the surface F or as a sum of isotopy classes in
the surface F mod orientation preserving homeomorphisms of F.)  We may
also compute the surface bracket polynomial by applying the skein
relation with the axiom: $ \langle F, U \rangle = d =-A^2 -A^{-2} $ if
$ U $ bounds a disk in $ F $.  In \cite{vom}, Manturov introduces
related polynomial invariants of knots in 2-surfaces. The surface
bracket polynomial is invariant under the Reidemeister moves II and
III.

\begin{rem}After collecting coefficients of equivalent (homologous or isotopic) 
surface-state pairs there may be fewer than $ 2^n $ surface-state
pairs with non-zero coefficients.
\end{rem}

We focus on using the homology classes of the curves in the
surface-state pairs to determine if a virtual knot diagram is
non-trivial.  Since $ F $ is a closed, orientable surface, we can
write $ F = T_1 \sharp T_2 ... \sharp T_n $ where each $ T_i $ is a
torus.  The homology group, $ H_1 (F) $, is generated by $ \lbrace
[m_1] \ldots [m_n],[l_1] \ldots [l_n] \rbrace $ where $ [m_i] $ and $
[l_i] $ represent the homology class of the meridian and longitude of
the torus $ T_i $ respectively.  If $ \gamma $ is a curve in the
surface-state pair $ (F,s ) $ then either $ \gamma $ is homologically
trivial in $ F $ or $ \gamma $ is homologous to a curve of the form:
\begin{equation*}
  \underset{i=1 \ldots n}{ \sum} a_i [m_i] + b_i [l_i]
\end{equation*}
where $ a_i $ and $ b_i $ are relatively prime.

The relationship between a cancellation curve, surface-knot pair
 and surface-state
 pairs is expressed in the following lemma.
\begin{lem} \label{lem1}Let $ C $ be
 a cancellation curve for a representation, $ ( F, K ) $, 
of a virtual knot diagram, $ \hat{K} $. Then 
$ C $ is a cancellation curve for every
surface-state pair.
\end{lem} 

\textbf{Proof}\qua Suppose that $ C $ is a cancellation curve for $
(F,K ) $, but $ C $ is not a cancellation curve for some surface-state
pair $ (F,s) $ obtained from the $ (F,K) $.  This indicates that
either $ C $ intersects a curve in the surface-state pair $ (F,s) $ or
$ C $ bounds a disk in the surface $ F $. We assumed that $ C $ was a
cancellation curve for the representation, so $ C $ does not bound a
disk in $ F $.  If the curve $ C $ intersects the state $ s $ then
then the curve $ C $ also intersects the original diagram $ K $ as a
result of the definition in Figure \ref{fig:smooth}. Hence, $ C $ is
not a cancellation curve for this representation $ (F,K) $,
contradicting our original assumption. \endproof

We define the \emph{intersection number} of two oriented curves, $
\alpha $ and $ \beta $, in a surface $ F $ to be the intersection
number between the elements $ [ \alpha ] $ and $ [ \beta ] $ of the
homology classes $ H_1 ( F, \mathbb{Z} ) $. We will denote this as $ [
\alpha ] \bullet [ \beta ] $. Recall from \cite{bredon} that
intersection number is the Poincare dual to the cup product, and that
it can be calculated by placing the two curves transversely to each
other and counting the sum of the oriented intersections of them.

In the next theorem we give conditions for the existence of a
cancellation curve $C$ for a representation $(F, K)$ of a virtual knot
diagram $\hat{K}$.  The conditions are given in terms of intersection
numbers of the corresponding surface state pairs with generators of
homology of the surface.  By applying these conditions to the surface
state pairs having non-zero coefficients in the surface bracket
polynomial, we find conditions for the existence of cancellation
curves for all diagrams obtained from $(F, K)$ by isotopy and
Reidemeister moves, in other words conditions for the minimality of
$(F,K)$.

\begin{thm}Let $ (F, K ) $ be a representation of a virtual knot diagram with 
$ F = T_1 \sharp T_2 \ldots 
\sharp T_n $. 
Let 
$$ \lbrace ( F, s_1), (F, s_2) \ldots (F,s_m) \rbrace $$ 
denote  the 
collection of surface-state pairs obtained from $ (F,K) $. 
Assign an arbitrary orientation to each curve in the surface-state
pairs.
Let $p: F \rightarrow T_k $ be the collapsing map, and let
$ p_*: H_1 (F, \mathbb{Z}) \rightarrow H_1 (T_k, \mathbb{Z} ) $
be the induced map on homology.  
If for each $ T_k $ there exist two states $ s_i $ and 
$ s_j $  with non-zero coefficients that contain curves 
(with arbitrarily assigned orientation) $ \gamma_i $ and $ \gamma_j $ respectively, 
such that $ p_*[ \gamma_i] \bullet p_*[ \gamma_j]  \neq 0 $ then
there is no cancellation curve for $ (F,K) $.
\end{thm}

\textbf{Proof}\qua We initially  assume that $ F $ is a torus, $ T $.
Note that $ p: F \rightarrow T $ is the identity map in this case.
Suppose $ (T,K) $ has a cancellation curve
 $ C $. Let $s_i $ and $ s_j $ 
be two states with non-zero coefficients 
that contain curves with non-zero intersection number in the torus $ T $.
The curve $ C $ is a cancellation curve and therefore does not
 intersect  any curve in the two states, $ s_i $ or $ s_j $. 
In the torus, each state consists of
 a collection of curves that are  parallel copies of a simple closed curve with non-trivial homology class in $ H_1 (T, \mathbb{Z} ) $ and curves that bound a disk in the surface $ T $ after homotopy.
If we arbitrarily assign an  orientation 
to the non-trivial curves, the
curves are either elements of the same homology class or cobound an 
annulus. 
Let $ [m] $ and $ [l] $ represent the homology classes 
containing the  oriented meridian and longitude respectively so that
 $ [m ] $ and $ [l] $ generate $  H_1 (T, \mathbb{Z} ) $.
Recall that if $ [ \gamma ] = a [m] + b [l] $ is a simple closed curve 
in $ T $ then $ a $ and $ b $ are relatively prime.

Let the state $ s_i $ contain $ \gamma_i $, a simple closed curve that does 
not bound a disk. Then the  homology class of  $ [ \gamma_i ]  $
is given by the equation $ [ \gamma_i ] 
= a_i [m] + b_i [l] $, where $ a_i $ and $ b_i $ are relatively 
prime. Let also the state 
$ s_j $ contain $ \gamma_j $, a simple closed curve that does not 
bound a disk such that
$ [ \gamma_j ] = a_j [m] + b_j [l] $.
By hypothesis, $  [ \gamma_i ] \bullet [ \gamma_j]
\neq 0 $, implying that the curves $ \gamma_i $ and $ \gamma_j $ 
are not elements of the same 
cohomology class and the curves do not cobound an 
annulus. Using homology classes, we compute that
$$
 [ \gamma_i ] \bullet [ \gamma_j ] = a_i b_j -  b_i a_j.
$$
The cancellation curve $ C $ is a simple closed curve that does 
 not bound a disk in $ T $. Let 
 $ [C]  = g [m] + f [l] $ where $ g $ and $ f $ are relatively prime.  
We note that by Lemma \ref{lem1}, the curve $ C $ does not intersect the 
curve $ \gamma_i $  or the curve 
$ \gamma_j $ since $ C $ is a cancellation curve. 
We compute that 
\begin{equation} \label{cgb1}
 [ \gamma_i ] \bullet [ C ] = a_if - g b_i = 0 
\end{equation}
and
\begin{equation} \label{cgb2}
 [ \gamma_j ] \bullet [ C ] = a_j f - g b_j = 0.
\end{equation}
Note that
$$
0= a_if - g b_i = a_j f - g b_j 
$$
and so
\begin{equation} \label{pf1}
 f(a_i - a_j) + g(b_j - b_i) = 0. 
\end{equation}
We will consider the following three possibilities: 
$ f \neq 0 $ and $ g \neq 0 $, $ f=0 $, or $ g=0 $. 

If we assume that $ f=0 $ then from \ref{cgb1} and \ref{cgb2} we obtain:
\begin{alignat*}{2} 
 -gb_i &= 0 & \qquad -gb_j &= 0 
\end{alignat*}
and as a result $ g = 0 $ or $ b_i = b_j = 0 $.
If $ g= 0 $ then $ C $ is not a cancellation curve 
because $ C $ bounds a disk in the torus. If $ b_i = b_j = 0 $ then
this   contradicts  the fact that
\begin{equation*}
  [ \gamma_i ] \bullet [ \gamma_j] \neq 0.
\end{equation*} 
Thus $ f = 0 $ is not possible.

By the same argument, $ g=0 $ is not possible.

Suppose that $ f \neq 0 $ and $ g \neq 0 $.
Recall that $ g $ and $ f $ are relatively prime, so that 
if $ f \neq 0 $ and $ g \neq 0 $ then either $ f = 1 $ or
 $ \frac{g}{f} $ is not an element
of the integers. From \ref{cgb1} and \ref{cgb2} we obtain:
\begin{alignat*}{2}
 a_i &=  g \frac{b_i}{f} & \qquad a_j &= g \frac{b_j}{f}
\end{alignat*}
Note that $ g $, $ a_i $, and $ a_j $ are integers and
the pairs $ (g,f) $, $ (a_i,b_i) $ and $ (a_j,b_j) $ are relatively
prime. This implies
that $ \frac{b_i}{f} $ and $ \frac{b_j}{f} $ are integers.
However, if $ \frac{b_i}{f} $ is an integer $ w $ such that 
$ w \neq \pm 1 $ then $ a_i = gw $ and $ b_i = fw $. This
contradicts the fact that the pair $ (a_i, b_i ) $ was relatively prime. We obtain a similar result for the pair $ ( a_j, b_j ) $. As a result, we determine that
\begin{equation*} 
\pm 1 = \frac{b_i}{f} = \frac{b_j}{f}.
\end{equation*}
Hence, 
$ b_i =  \pm f $ and $ b_j = \pm f $. 
Correspondingly, $ a_i = \pm g $ and $ a_j = \pm g $. 
This contradicts the fact that $ [ \gamma_i ] \bullet [ \gamma_j] \neq 0 $.
 
Therefore, $ C $ is not a cancellation curve for the torus $ T $.

Let $ F= T_1 \sharp T_2 \sharp \ldots \sharp T_n $. Let C be a 
cancellation curve for the surface $ F $. Let
 $ [m_k] $ and $ [l_k] $ represent the homology classes containing 
the meridian and longitude of the torus $ T_k $ respectively. 
Let $ p_* H_1 (F, \mathbb{Z} ) \rightarrow H_1 ( T_k, \mathbb{Z}) $.
We note that $ p_*( [C] ) = f [m_k] + g [l_k] $ with either $ f \neq 0 $ or
$ g \neq 0 $ for some $ T_k $. Otherwise $ p_*[ C ] $ would bound a 
disk in each 
$ T_k $.  As a result, the curve $ C $ divides the surface $ F $ into 
two components, one of which contains the knot.
If $ C $ bounds a disk, then $ C $ is not a cancellation curve. Hence 
$ C $ bounds a component containing some states $ s_i $ and $ s_j $, 
contradicting the fact that $ C $ is a cancellation curve.

Let $ s_i $ and $ s_j $ be states such that $ p_* [s_i] $ contains 
a curve $ \gamma_i $ and $ p_* [s_j] $ contains a curve $ \gamma_j $
such $ [ \gamma_i ] \bullet [ \gamma_j] \neq 0 $.
Let $ [ \gamma_i] = a_i [m_k] + b_i [l_k] $ and let 
$ [ \gamma_j ] =  a_j [m_k] + b_j [l_k] $ where the pairs
$ (a_i, b_i) $ and $ ( a_j, b_j ) $ are relatively prime.
Using the argument given previously, we eliminate the possibility
that
$ f = 0 $ or $ g = 0 $. We then consider the 
cases when $ f \neq 0 $ and $ g \neq 0 $.
Using the same argument, we determine that
$$
 b_i =\pm f \text{ or } b_j = \pm f
$$
and combined with \ref{pf1} this
 indicates that $ a_i = \pm g $ and $ a_j = \pm g $, contradicting our 
assumption that $ [ \gamma_i ] \bullet [ \gamma_j] \neq 0  $.
The cancellation curve $ C $  was arbitrary and therefore $ F $ has no
cancellation curves.\endproof

\begin{rem} We note that the condition of theorem corresponds to the existence
of two non-trivial, non-isotopic curves in each torus component projected from the states of the representation $ (F, K ) $.
\end{rem}

\section{Virtual knot diagrams with one virtualized 
crossing}

Recall that a representation of a virtual knot diagram is minimal if
no handles can be removed after a sequence of Reidemeister moves.  In
this section, we use the surface bracket polynomial to prove
minimality for a class of virtual diagrams with one virtualized
crossing. This enables us to show that many virtual knot diagrams are
non-classical.

\begin{figure}[htb] \epsfysize = 1.3 in
\centerline{\epsffile{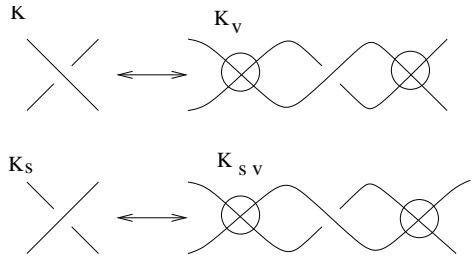}}
\caption{Virtualized Crossing}
\label{fig:virtualized}
\end{figure}

A classical crossing in a virtual knot diagram is \emph{virtualized} by the 
following procedure: a tangle consisting of a single crossing is removed and 
replaced with a tangle consisting of the opposite crossing flanked by two 
virtual crossings. This procedure is illustrated in Figure \ref{fig:virtualized}.

\begin{figure}[htb] \epsfysize = 1 in
\centerline{\epsffile{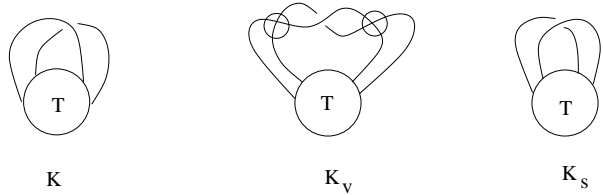}}
\caption{$K$, $ K_s $ and the virtualized diagram: $K_v $}
\label{fig:virted}
\end{figure}

We consider the three knot diagrams as shown in Figure  \ref{fig:virted}. The
first diagram, labeled $ K $ is a classical knot diagram formed by 
one isolated classical crossing $ v $ and the classical tangle $ T $.
The second diagram, $ K_v $ is obtained from $ K $ by virtualizing 
the crossing $ v $. The third diagram, $ K_s $ is obtained by switching the isolated crossing.

\begin{figure}[htb] \epsfysize = 0.3 in
\centerline{\epsffile{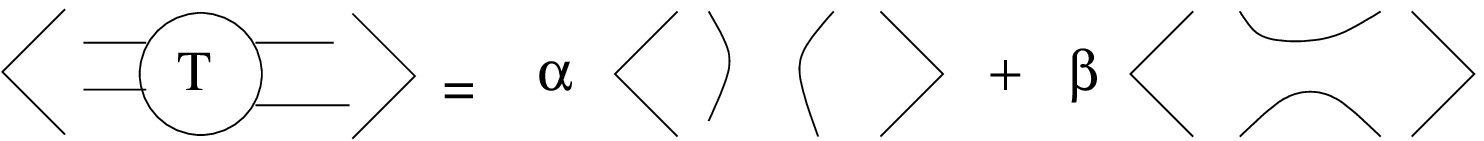}}
\caption{Skein Relation of Tangle T }
\label{fig:skeint}
\end{figure}

We apply the skein relation to the tangle $ T $ and obtain the 
relation shown in Figure \ref{fig:skeint} where $ \alpha $ and 
$ \beta $ 
are coefficients in $ \mathbb{Z} [A, A^{-1} ] $.

Applying the bracket skein relation and the relation shown in Figure
\ref{fig:skeint}, we determine that:
\begin{align} \label{skeink}
\langle K \rangle &= -A^{-3} \alpha - A^3 \beta \\
\langle K_s \rangle &= -A^{3} \alpha  - A^{-3} \beta. \nonumber
\end{align}
We note that $ \langle K_v \rangle = \langle K_s \rangle $, see \cite{kvirt}.
Hence
\begin{align*} \label{skeinkv}
\langle K_v \rangle &= - A^{3} \alpha  - A^{-3} \beta. 
\end{align*}
In particular, if $ K_s $ is an unknot and $ w $ is the writhe of $ K_s $. Then
$ \langle K_s \rangle = (-A)^{-3w} $ and $V_{K_v}(t)= V_{K_s} (t) = 1 $. 
 
\begin{lem}Let $ K $ be a classical knot or link, and let $ K_s $ and 
$ K_v $ be as in Figure \ref{fig:virted}. Let 
$ \alpha $ and $ \beta $ be defined as in Figure \ref{fig:skeint}. Then
\begin{align*}
\langle K \rangle &= A^{-6} \langle K_s \rangle
  +( -A^3 + A^{-9}) \beta \\
\alpha &= A^{-3} \langle K_s \rangle  + A^{-6} \beta
\end{align*}
\end{lem}
\textbf{Proof}\qua 
Using the second part of equation \ref{skeink},
\begin{equation*} 
  \langle K_s \rangle = -A^{3} \alpha  -A^{-3} \beta. 
\end{equation*}
Solving for $ \alpha $, we determine that:
\begin{equation*}
   \alpha = -A^{-3} \langle K_s \rangle - A^{-6} \beta. 
\end{equation*}
Substitute into equation \ref{skeink} and find:
\begin{align*}
 \langle K \rangle &= A^{-6} \langle K_s \rangle + 
	( -A^3 + A^{-9}) \beta \tag*{\qed}
\end{align*}
We introduce the following proposition.
\begin{prop} \label{zerocor}  
Let $ K $ be a classical knot or link, and let $ K_s $ and 
$ K_v $ be as in Figure \ref{fig:virted}. Let 
$ \alpha $ and $ \beta $ be defined as in Figure \ref{fig:skeint}. Then
 $ \langle K \rangle = ((-A)^{3})^{ \pm 2} \langle K_s \rangle $ 
  if and only if $ \alpha =0 $ or $ \beta =0 $.
\end{prop}
\textbf{Proof}\qua
Suppose that $ \langle K \rangle = ((-A)^3)^{ \pm 2} \langle K_s \rangle $.
We compute that:
\begin{gather*}
\langle K \rangle  = -A^{-3} \alpha - A^{3} \beta \\
\text{ and } \\
 \langle K_s \rangle = -A^{3} \alpha - A^{-3} \beta
\end{gather*}
where $ \alpha $ and $ \beta $ are non-zero elements of 
$ \mathbb{Z}[A,A^{-1}] $ as shown in Figure \ref{fig:skeint}.
As a result, we observe that: 
\begin{gather*}
    \langle K \rangle =  - A^{-3} \alpha - A^{3} \beta \\
\text{  and   } \\
  A^{ \pm 6 } \langle K \rangle = - A^{3} \alpha - A^{-3} \beta  
\end{gather*}  
Now, taking $ + 6 $ and $ - 6 $ respectively, we find:
\begin{gather*}
  \langle K  \rangle = -A^{-3} \alpha - A ^{3} \beta \text{ and } \\
  \langle K \rangle = -A^{-3} \alpha - A^{-9} \beta \\
\text{or} \\
\langle K  \rangle = -A^{-3} \alpha -A^{3} \beta   \text{ and } \\
  \langle K \rangle = -A^{9} \alpha - A^{3} \beta 
\end{gather*}
These equations are contradictory unless either $ \alpha = 0 $ or $ \beta =0 $.  
Suppose that $ \alpha = 0 $. Using the skein relation, 
\begin{gather*}
\langle K \rangle = -A^3 \beta \\
\text{ and } 
\langle K_s \rangle = -A^{-3} \beta
\end{gather*}
Therefore, $ \langle K \rangle = A^6 \langle K_s \rangle $. 
We may perform a similar computation if $ \beta = 0 $ and determine that 
$ \langle K \rangle = A^{-6} \langle K_s \rangle $.
\endproof

Note that this proposition tells us that if $ K_s $ is an unknot or an unlink then $ K $ has the same bracket polynomial as an unknot or unlink if and only if $ \alpha =0 $ or $ \beta =0 $.

\begin{figure}[htb] \epsfysize = 1 in
\centerline{\epsffile{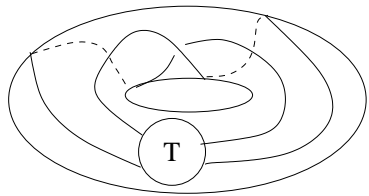}}
\caption{Representation: $(F, K_v) $}
\label{fig:torusrep}
\end{figure}

We consider a representation of the virtual knot diagram $ K_v $ as a knot or a link
embedded in a torus $ F $  shown in 
Figure  \ref{fig:torusrep}.

\begin{thm} \label{tst} 
Let $ K $ be a classical knot or link diagram as in Figure \ref{fig:virted}
with associated links $ K_s $ and $ K_v $. If $ \alpha $ and $ \beta $, as determined in Figure \ref{fig:skeint}, are both non-zero then $ K_v $ is a non-classical and non-trivial virtual link.
\end{thm}

\begin{figure}[htb] \epsfysize = 1 in
\centerline{\epsffile{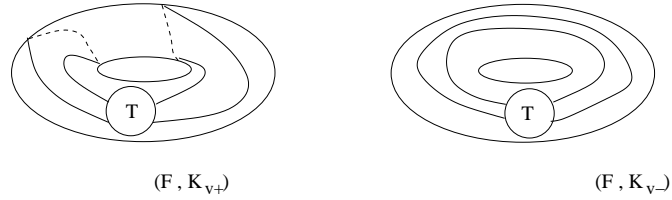}}
\caption{States in the Torus}
\label{fig:torusstates}
\end{figure}

\textbf{Proof}\qua We obtain the  two
surface-state pairs $(F, K_{v+} ) $ and $ (F,K_{v-} ) $
 in Figure \ref{fig:torusstates} from the skein relation.

\begin{figure}[htb] \epsfysize = 2.5 in
\centerline{\epsffile{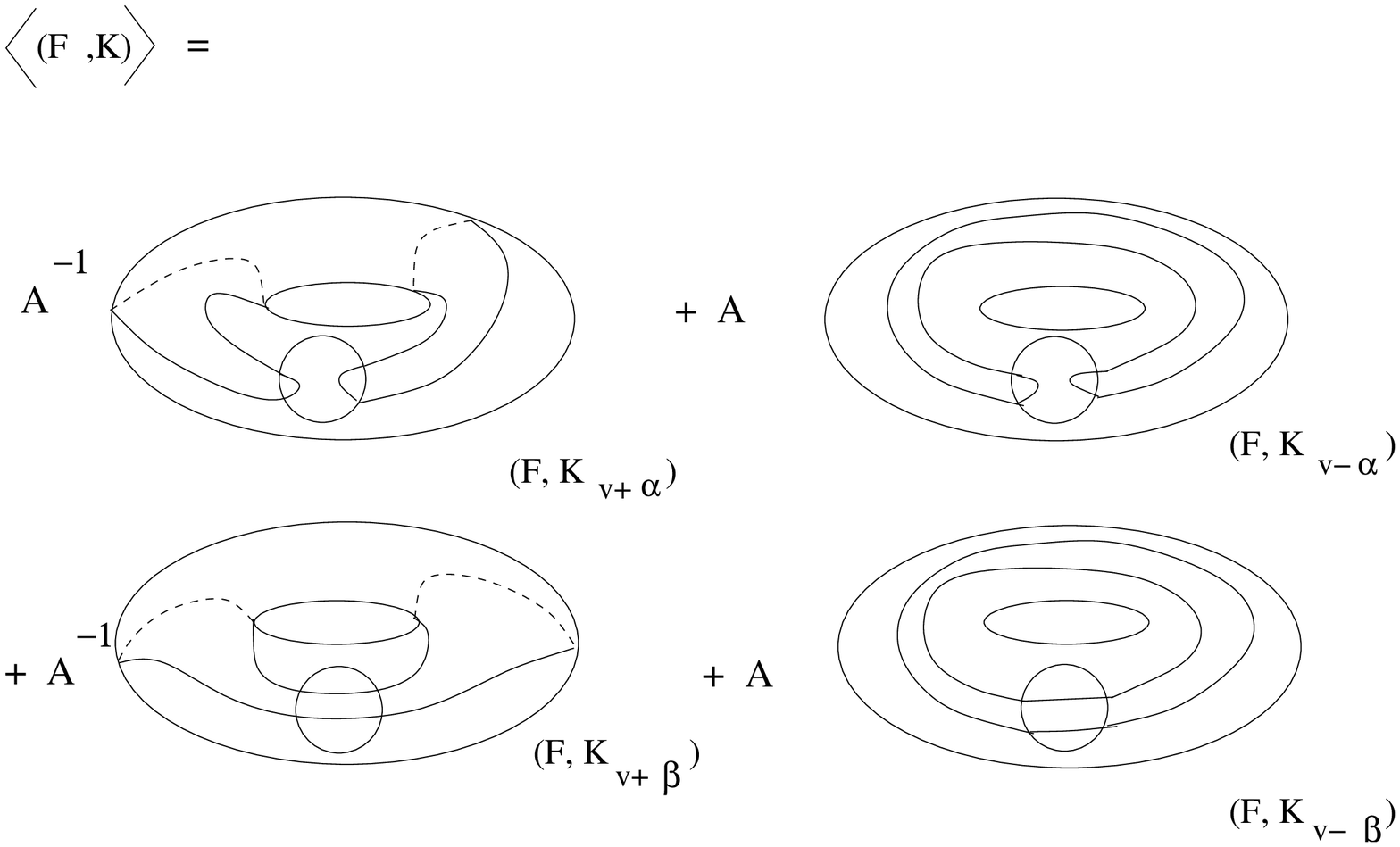}}
\caption{Surface-State Equation}
\label{fig:kcalc}
\end{figure}

Hence, 
\begin{equation*}
\langle (F,K_v) \rangle = A^{-1}  \langle (F,K_{v+}) \rangle
 + A \langle (F, K_{v-} ) \rangle 
\end{equation*}
Combining this expansion with that states from \ref{fig:torusstates},
we obtain the relation shown in Figure \ref{fig:kcalc}.

We note that 
\begin{gather} \label{tstates}
\langle (F,K_v) \rangle = A^{-1}( \alpha  \langle (F, K_{v+,\alpha})
 \rangle + \beta \langle (F, K_{v+, \beta} ) \rangle ) \\
+ A ( 
 \alpha  \langle (F, K_{v-,\alpha})
 \rangle + \beta \langle (F, K_{v-, \beta} ) \rangle ). \nonumber
\end{gather}
Referring to Figure \ref{fig:kcalc}, we observe that the states 
$ (F,K_{v-, \alpha} ) $ and $ ( F, K_{v+, \beta} ) $ both contain a single 
curve that bounds a disk in $ F $. As a result equation \ref{tstates} 
reduces to
\begin{equation} \label{tstates2}
\langle (F, K) \rangle = (A \alpha + A^{-1} \beta) + (A^{-1}  \alpha \langle 
(F, K_{v+, \alpha} ) \rangle + A \beta \langle (F, K_{v-, \beta}) \rangle.
\end{equation}
Note that if both
 $ \alpha $ and $ \beta $ are non-zero, the subspace 
of curves generated by 
the surface-states spans the space of curves in the torus. \endproof

 We recall the following theorem from \cite{heckalg}:
\begin{thm}[V.F.R. Jones] \label{vjones} If $ K $ is a knot then $ 1- V_K (t) = W_K (t) (1-t) ( 1-t^3)$ for some Laurent polynomial $ W_K(t) $.
\end{thm}
Note that if $ V_K (t) $ is a monomial then $ V_K (t) = 1 $. If $ V_K (t) $ is a monomial, $(V_K(t)= t^n )$, then $ 1 - t^n $ is divisible by $ (1-t) $ and $ (1-t^3)$ by Theorem \ref{vjones}. 
Hence, if K is a knot diagram and $ \langle K \rangle = (-A)^n $ then 
$ n = -3w $, where $ w $ is the writhe of $ K $. 
We obtain the following corollary from this fact.

\begin{cor} \label{nontriv} If $ K $ is a classical knot diagram with unknotting number one and non-unit Jones polynomial and  
 $ K_s $ is the unknot then $ K_v $ is non-classical and non-trivial.
\end{cor}
\textbf{Proof}\qua 
Let  $ K $ have writhe $ w  $ then $ K_s $ is the 
unknot with writhe $ w \pm 2 $.  We obtain:
$ \langle K_S \rangle =  (-A)^{-3(w \pm 2)} $.
 By Corollary \ref{zerocor} 
 $\alpha =0 $ or $ \beta =0 $ 
if and only if $ \langle K \rangle = (-A)^{-3w} $.
Since $ \langle K \rangle \neq (-A)^{-3w} $ then $ \alpha \neq 0 $ and $ \beta \neq 0 $. 
This indicates that the given representation of $ K_v $ has no cancellation curves.
The virtual genus of $ K_v $ is one, indicating that $ K_v $ is non-classical 
and non-trivial. \endproof

\begin{rem} Note that Corollary \ref{nontriv} does not eliminate the possibility that 
there exists a non-trivial classical knot diagram $ K $
where both $ K $ and $ K_s $  have  unit Jones polynomial,  but $ K_v $ is not detected by the surface bracket 
polynomial.
\end{rem}

In \cite{silwill}, the following theorem is obtained from an 
analysis of the fundamental group.

\begin{thm}[Silver--Williams] \label{silwilth} Let $ K $ be a non-trivial classical knot diagram, and $ v $ is 
a classical crossing. If $ K_v$ is the virtual knot diagram 
obtained by virtualizing $ v $ in $ K $ then $ K_v $ is non-classical 
and non-trivial. \end{thm}

 If $ K $ is a non-trivial classical knot with $ V_K (t) = 1 $ and   $ V_{K_s} (t) = 1 $ then the surface bracket 
polynomial would not detect $ K_v $ even though the virtual genus is one via Theorem \ref{silwilth}. 

\begin{figure}[htb] \epsfysize = 1.4 in
\centerline{\epsffile{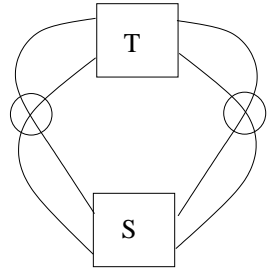}}
\caption{A virtual knot diagram
constructed from two tangles}
\label{fig:twotangle}
\end{figure}

We may generalize our procedure to demonstrate that a larger class of virtual knot diagrams is non-trivial and non-classical.  
Construct a virtual knot diagram from two classical tangles, $ T $ and $ S $ as shown in Figure \ref{fig:twotangle}.
The same arguments prove that if the tangles are expanded as shown in Figure \ref{fig:skeint} and the coefficients, $ \alpha $ and $ \beta $, are non-zero for both tangles then the virtual knot diagram has virtual genus one, \textit{whence the virtual knot diagram is non-classical and non-trivial.}

\begin{figure}[htb] \epsfysize = 1.5 in
\centerline{\epsffile{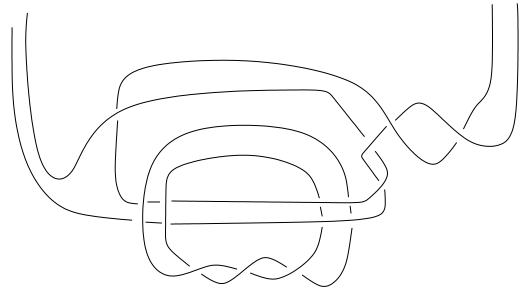}}
\caption{A tangle $ T_L $ resulting in a link}
\label{fig:linktangleT}
\end{figure}

\begin{figure}[htb] \epsfysize = 3 in
\centerline{\epsffile{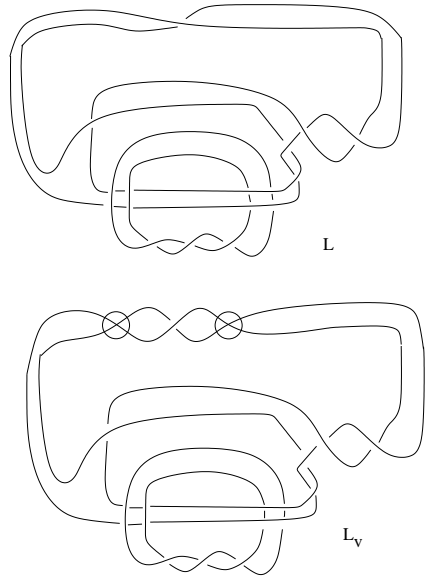}}
\caption{Link Diagrams: $ L $ and $ L_v $}
\label{fig:linkseries}
\end{figure}

\subsection*{Undetectable Examples}
Given the tangle shown in Figure \ref{fig:linktangleT}, we use the method given at the beginning of this section to construct a link diagram by taking a tangle sum with a single crossing.
The link constructed from this tangle is shown in Figure \ref{fig:linkseries}. This link, $ L $  
and the corresponding link $L_s$ with a switched crossing  have the property that both $ L $ and $ L_s $ have the same Jones polynomial as an unlink of two components. These link diagrams were constructed using the methods of \cite{morwen}.

For this link $ L $, 
we note that $ L_v $ is not detected by the surface bracket polynomial since by Corollary
\ref{zerocor}, $ \alpha = 0 $. We thank Alexander Stoimenow for pointing out the usefulness of \cite{morwen}.

\begin{figure}[htb] \epsfysize = 1 in
\centerline{\epsffile{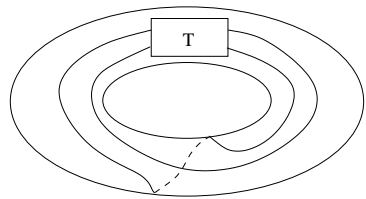}}
\caption{Representation of virtual knot diagram with 1 virtual crossing}
\label{fig:sampletorus}
\end{figure}

\begin{rem}
We briefly comment on the case of a virtual knot diagram with a single virtual crossing and a classical tangle $ T $. Let $ K $ be such a virtual knot diagram. A schematic representation of this virtual knot diagram $ (F,K) $ is shown in Figure \ref{fig:sampletorus}.
Let $ [m] $ and $ [l] $ represent the meridian and longitude of the torus $ F $. 
If we expand the tangle T as illustrated in Figure \ref{fig:skeint}, we obtain:
\begin{gather*}
\langle K \rangle = \alpha + \beta \\
\text{and} \\
\langle (F,K) \rangle = \alpha \langle (F,[m]) \rangle + \beta\langle
( F, [m+2l]) \rangle
\end{gather*}
Note that if $ \alpha \neq 0 $ and $ \beta \neq 0 $ then $ K $ is non-classical and non-trivial. If $ \alpha = 0 $ or $ \beta =0 $ then no decision can be made. In particular, we can construct a virtual link diagram with a single virtual crossing using the tangle shown in Figure \ref{fig:linktangleT}. This link is not detected by the surface bracket polynomial. 
\end{rem}

\section{Other virtual knot diagrams}

We study other virtual knot diagrams and determine if these diagrams are 
non-classical
and hence non-trivial using this technique.

\begin{figure}[htb] \epsfysize = 1 in
\centerline{\epsffile{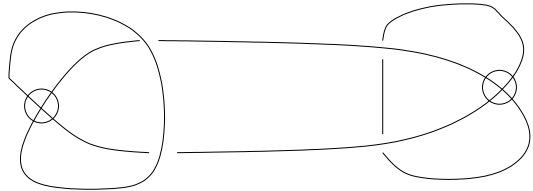}}
\caption{Kishino's Knot }
\label{fig:kishknot}
\end{figure}

Kishino's knot is illustrated in Figure \ref{fig:kishknot}.

\begin{figure}[htb] \epsfysize = 1 in
\centerline{\epsffile{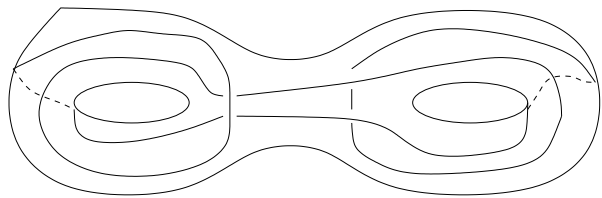}}
\caption{Genus Two Representation  of Kishino's Knot }
\label{fig:charkish}
\end{figure}

This knot has a trivial fundamental group and bracket polynomial. 
 This knot was determined to be non-trivial by the 3-strand bracket polynomial 
\cite{kishpoly} 
and the quaternionic biquandle \cite{BIQ}. Both of these methods involve 
intensive and difficult computation.
The methods introduced in this paper demonstrate 
that Kishino's knot is non-classical and non-trivial and 
 that the virtual genus of Kishino's knot is greater than zero.
Recall the definition of virtual genus as given before Remark \ref{virtgen}.
\begin{thm}The virtual genus of Kishino's knot 
is two. 
\end{thm}
\begin{cor} Kishino's knot is non-trivial and non-classical. 
\end{cor}
\textbf{Proof}\qua
We show a genus two representation of Kishino's knot in 
Figure \ref{fig:charkish}.

Note that Kishino's knot has 4 crossings. By 
application of the bracket polynomial, we 
obtain 16 surface-states from this representation.

We illustrate the 16 surface states in the following Figures.
\begin{figure}[htb] \epsfysize = 2 in
\centerline{\epsffile{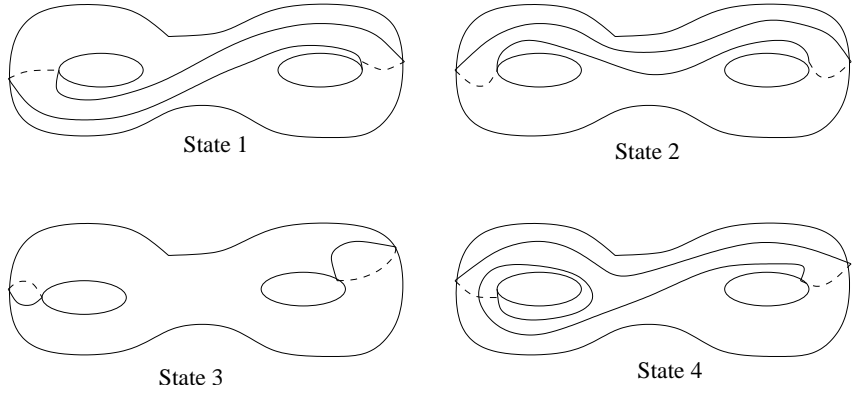}}
\caption{States of Kishino's Knot, 1-4 }
\label{fig:kishstates1}
\end{figure}
\begin{figure}[htb] \epsfysize = 2 in
\centerline{\epsffile{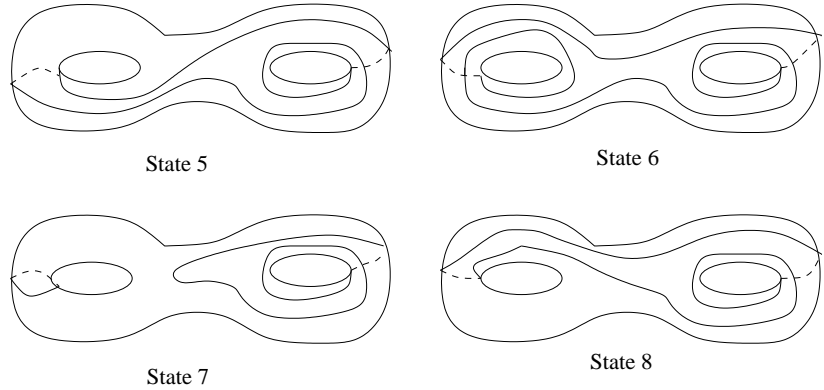}}
\caption{States of Kishino's Knot, 5-8 }
\label{fig:kishstates2}
\end{figure}
\begin{figure}[htb] \epsfysize = 2 in
\centerline{\epsffile{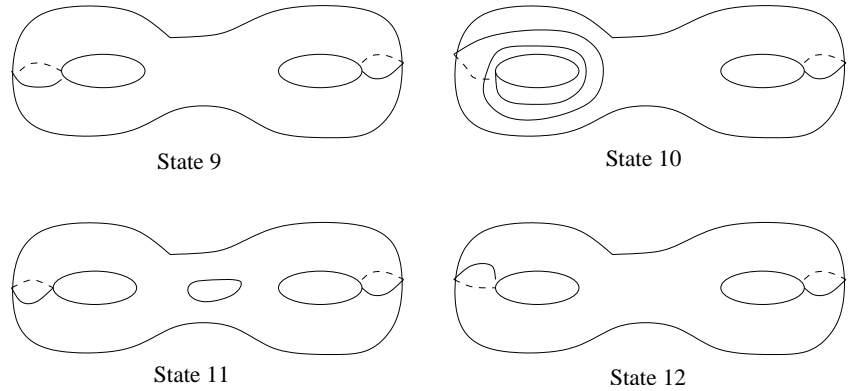}}
\caption{States of Kishino's Knot, 9-12 }
\label{fig:kishstates3}
\end{figure}
\begin{figure}[htb] \epsfysize = 2 in
\centerline{\epsffile{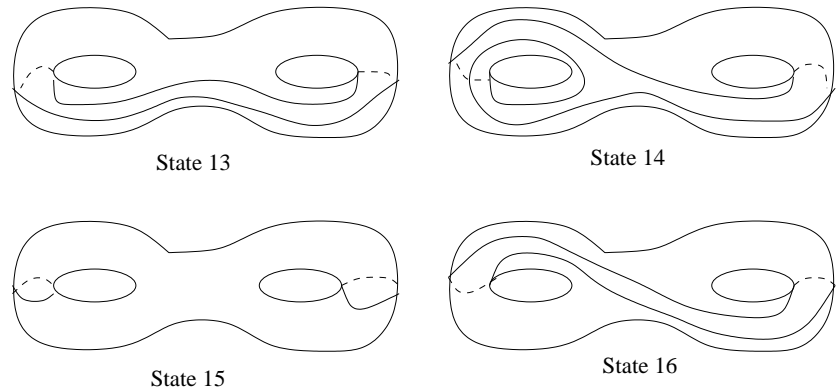}}
\caption{States of Kishino's Knot, 12-16 }
\label{fig:kishstates4}
\end{figure}

We will denote state i as $s_i $ and the coefficients as $c_i $.
\begin{alignat}{4}
c_1 &= 1 & \qquad c_2 &= A^4 & \qquad c_3 &= A^2 & \qquad c_4 &= A^2 \\
c_5 &= A^{-2} & \qquad c_6 &= 1 & \qquad c_7 &= 1	 & \qquad c_8 &= A^2 \\
c_9 &= A^{-2} & \qquad c_{10} &= 1 & \qquad c_{11} &= 1 
& \qquad c_{12} &= A^2 \\
c_{13} &= A^{-4} & \qquad c_{14} &= A^{-2} & \qquad c_{15} &= A^{-2} & \qquad c_{16} &= 1
\end{alignat}
We combine states with isotopy curves and obtain the following formula for the
surface bracket polynomial.
\begin{gather*}
  (F,s_1) + A^4(F,s_2)+ (A^2+A^{-2})(F,s_3) + A^2(F,s_4)\\
  + A^{-2} (F,s_5) + (F,s_6) + (F,s_7) + A^2(F,s_8)  \\
+ (F,s_{10}) + A^{-4}(F,s_{13}) + A^{-2}(F,s_{14}) + (F, s_{16})
\end{gather*}
Note that the states $ s_3, s_{10}, $ and $ s_{14} $ modulo 2 span the entire space of homology classes of curves in the connected sum of two tori. 
The fact that these curves span the homology group is invariant under isotopy of the knot in the surface and invariant under homeomorphisms of the surface. Hence, Kishino's knot is not equivalent to a knot that admits a cancellation curve.
Therefore, the virtual genus of Kishino's knot is two. \endproof

\begin{figure}[htb] \epsfysize = 1 in
\centerline{\epsffile{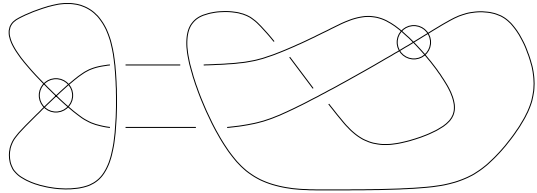}}
\caption{Modified Kishino's Knot}
\label{fig:halfkish}
\end{figure}

We consider a slight modification of Kishino's knot, as illustrated
in Figure \ref{fig:halfkish}.

This virtual knot diagram is undetected by the fundamental group and the
1-strand and 2-strand bracket polynomial \cite{dye}.
The knot has 6 classical crossings, and expands into 32 states. 
We may represent this virtual knot diagram as a knot diagram on 
 the connected sum of two tori, 
and compute the expanded states  and coefficients in each torus. 
This process forces
us to conclude that there are no cancellation curves in this surface. 
As a result we obtain:
\begin{prop}The modified Kishino's knot is non-trivial and non-classical. 
\end{prop}
\textbf{Proof}\qua
Use the method given in the previous proof.
 Compute the surface bracket polynomial and compare the rank of the
equivalence classes of curves in states with non-zero coefficients.\endproof

\begin{figure}[htb] \epsfysize = 1 in
\centerline{\epsffile{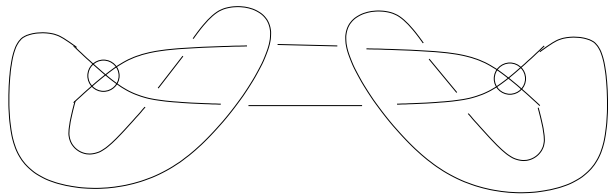}}
\caption{New Knot }
\label{fig:heather}
\end{figure}

We consider a further modification of this virtual knot diagram, as 
shown in Figure \ref{fig:heather}.

\begin{figure}[htb] \epsfysize = 1 in
\centerline{\epsffile{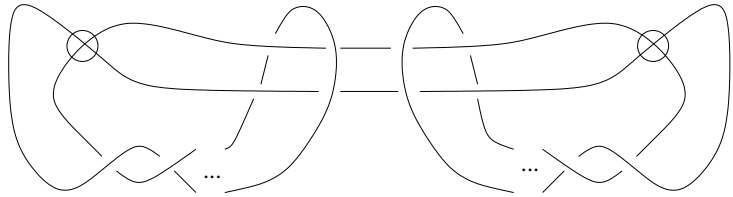}}
\caption{Schematic of the Family }
\label{fig:schem}
\end{figure}

\begin{thm} \label{modk} The virtual knot diagram shown in Figure 
\ref{fig:heather} is non-trivial 
and non-classical.
\end{thm}
In fact, this diagram is part of an infinite class of knots that is not detected by the bracket polynomial. 
We will prove that the members of the infinite class are detected by the surface bracket polynomial. This includes the case of Theorem
\ref{modk}.
\begin{thm} There is an infinite family of non-trivial virtual knot diagrams
obtained by modifying Kishino's knot. These virtual knot diagrams are not detected by the 
bracket polynomial but are detected by the surface
bracket polynomial. A schematic diagram of this family is shown
in Figure \ref{fig:schem}.
\end{thm}

\begin{figure}[htb] \epsfysize = 3 in
\centerline{\epsffile{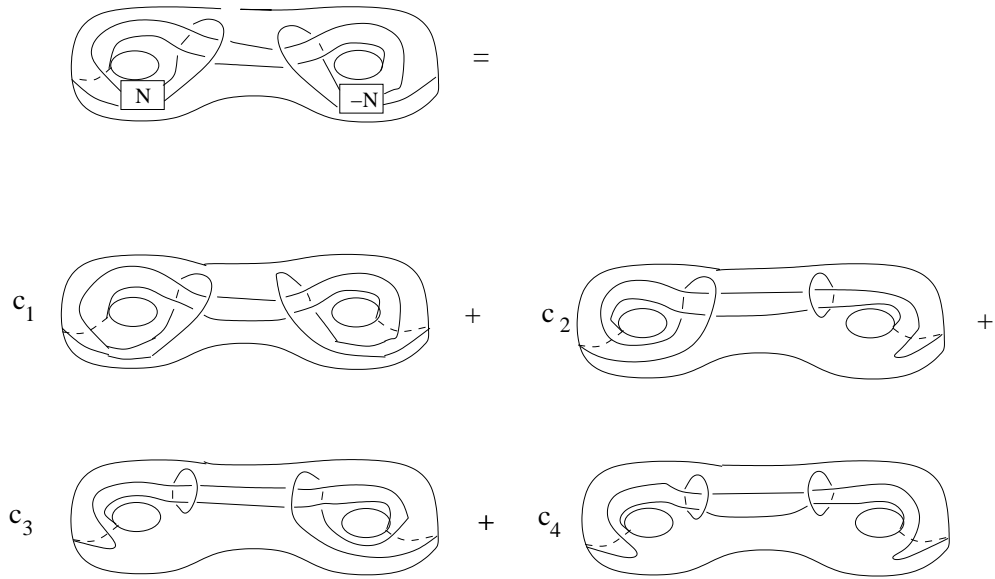}}
\caption{Partial State Expansion}
\label{fig:expan2}
\end{figure}

\textbf{Proof}\qua
We denote the members of this family as $ P_n $, where $ n $ denotes the number of inserted twists. As a result, $ P_0 $ refers to the diagram shown in Figure \ref{fig:heather}.
By applying the surface bracket polynomial to the knot shown in Figure \ref{fig:heather}, we obtain the following states with non-zero coefficients.
\begin{figure}[htb] \epsfysize = 2 in
\centerline{\epsffile{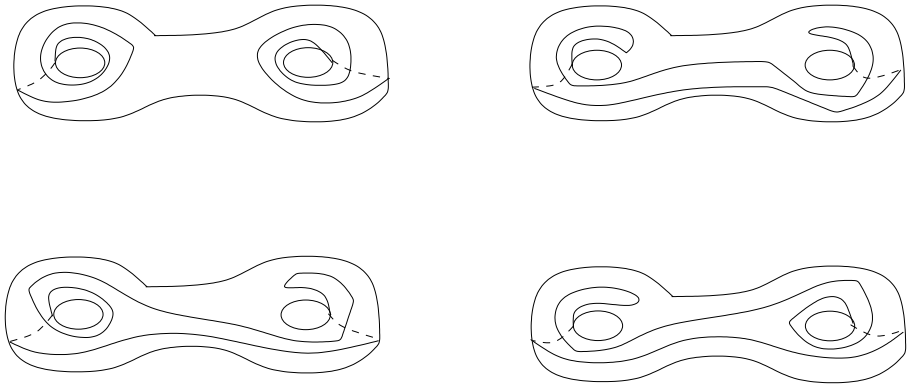}}
\caption{States of the New Knot shown in Figure \ref{fig:heather} }
\label{fig:expan1}
\end{figure}
These states are sufficient to ensure that no cancellation curves exist in the surface. Hence, the virtual genus of this diagram is two.
We expand the diagram $ P_n $   to obtain the state sum illustrated in Figure \ref{fig:expan2} using the skein. 
A lengthy calculation shows that the coefficients $ c_1, c_2, c_3 $ and $ c_4 $ are non-zero. The states shown in Figure \ref{fig:expan1} are obtained by expanding the state with coefficient $ c_1 $ from Figure \ref{fig:expan2}. The expansion of the states coefficients $ c_2, c_3 $ and $ c_4 $
does not involve these states. Consequently, the states shown in Figure \ref{fig:expan1} are not cancelled and have non-zero coefficients in the final state sum. These states are sufficient to ensure that no cancellation curves exist. Hence, these virtual knot diagrams have virtual genus two. As a result, they are non-trivial and non-classical. \endproof

The knot diagram in \ref{fig:heather} is not detected by the 1 and 2-strand bracket polynomial, but it is detected by the 3-strand bracket polynomial. However it is simpler 
to apply to bracket polynomial a genus 2 representation of this knot. 
We do not know if the other members of the infinite family are detected by the 3-strand bracket polynomial. These computations are extremely complex, and we are currently unable to complete the calculations on a computer. 

For a virtual knot diagram with $ n $ classical crossings, the 3-strand bracket polynomial 
has complexity of order $ 2^{9n} $.  These computations are considered in depth in \cite{dye}. 
\begin{rem}We conjecture that the 3-strand bracket polynomial
 detects virtual knot diagrams. The states of the 3-strand bracket 
polynomial may reflect the geometry of the minimal surface. 
\end{rem}
\begin{rem} Kodakami's work on the detection of virtual knot diagrams is closely 
related to this approach \cite{kodakami}. We note that his approach works 
for diagrams that are non-trivial in the flat category. The flat versions
of 
the virtual knot diagrams in Figure \ref{fig:schem} and in 
Figure \ref{fig:halfkish} are trivial, indicating that Kodakami's method would not detect these knots. 
\end{rem}

\section{Virtual knot diagrams with two virtualized 
crossings}

\begin{figure}[htb] \epsfysize = 2 in
\centerline{\epsffile{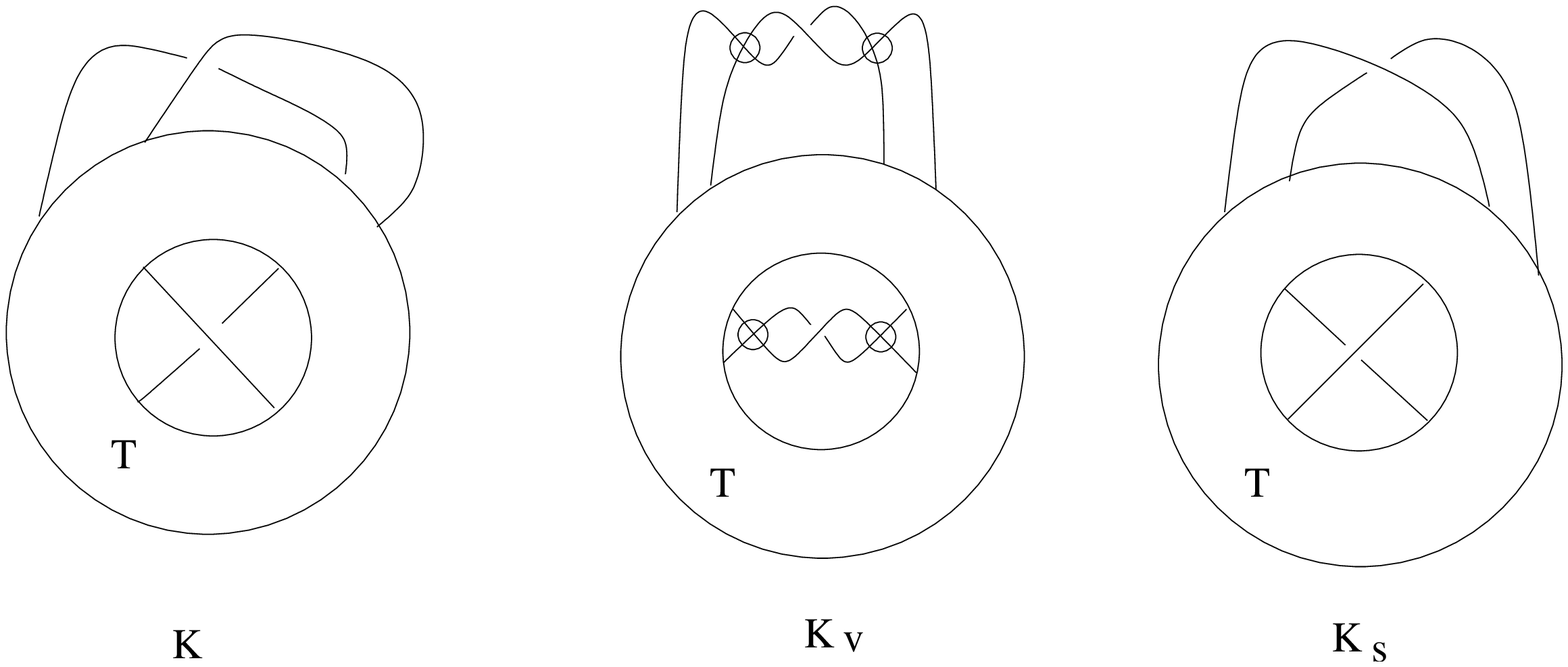}}
\caption{$ K $, $K_v$, and $ K_s $ }
\label{fig:virt2}
\end{figure}

We conclude this paper by considering the following class of virtual knot 
diagrams. Let $ K $ be a classical knot diagram, consisting of a
classical 4-4 tangle $ T $, occuring in an annulus, and two isolated crossings. The isolated crossings 
as chosen so the knot $ K_s $ with the isolated crossings switched in 
the unknot. Let $ K_v $ denote the modified diagram produced by 
virtualizing the isolated crossings. These Figures 
are illustrated in Figure \ref{fig:virt2}. Note that the genus of the characterization of $ K_v $ is bounded above 
by genus 2.

 We apply our new
method to a virtual knot diagram constructed by applying two virtualizations. Observe 
that in some cases it is possible to determine 
that a virtual knot diagram is non-classical and 
non-trivial without a full expansion of the bracket polynomial.

\begin{figure}[htb] \epsfysize = 2 in
\centerline{\epsffile{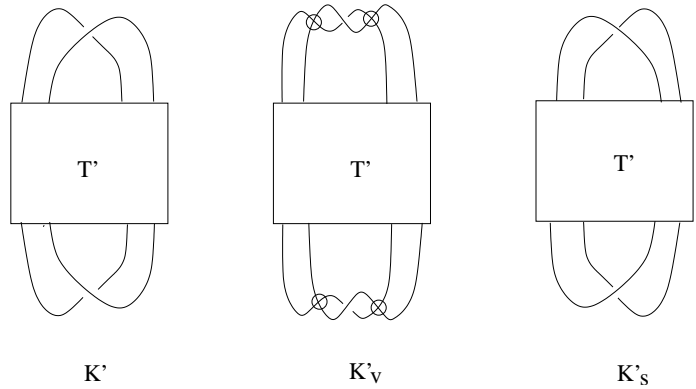}}
\caption{$K' $, $ K'_v $, and $ K'_s $ }
\label{fig:virted2}
\end{figure}

The equivalence classes of states that arise from application
of the bracket polynomial to
 4-4 tangle in an annulus have not been determined. As a result we 
restrict our attention to a modification of $ K $. The knot $ K' $ is 
obtained by applying a sequence of Reidemeister moves to the classical 
knot diagram $ K $. The diagram $ K' $ consists of a classical 4-4 tangle $ T' $, contained in a disk, and two isolated crossings. We construct $ K'_v $ and $ K'_s $ as before. The
diagrams $ K' $, $ K'_v $ and $ K'_s $ are illustrated in Figure 
\ref{fig:virted2}.

The diagrams $ K $ and $ K' $ are equivalent, but 
two virtual diagrams $ K_v $ and $ K'_v $ are not 
necessarily virtually equivalent. The diagrams $ K_s $ and $ K'_s$ are 
both unknots. 
We consider the bracket polynomial:
\begin{align*}
 \langle K \rangle &= (-A)^{3n} \langle K' \rangle \\
 \langle K_s \rangle &= (-A)^{3n} \langle K'_s \rangle \\
 \langle K_s \rangle &= \langle K_v \rangle \\
 \langle K'_s \rangle &= \langle K'_v \rangle
\end{align*}
where $ n $ reflects the number of Reidemeister I moves. 

\begin{figure}[htb] \epsfysize = 2.5 in
\centerline{\epsffile{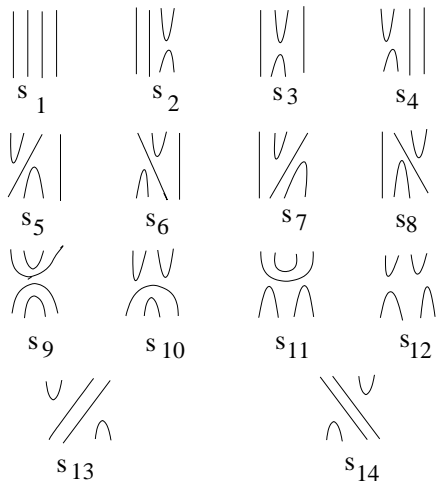}}
\caption{Generators of the $ 4^{th}$ Temperly-Lieb Algebra }
\label{fig:tl4}
\end{figure}

Expanding the tangle $ T' $ using the skein relation, we obtain a linear
combination of the twelve  
elements of the $ 4^{th} $ Temperly-Lieb algebra \cite{tl}
with coefficients in $ \mathbb{Z} [A, A^{-1}] $.
The twelve elements of the $ 4^{th} $ Temperly-Lieb algebra are shown 
in Figure \ref{fig:tl4}.

\begin{figure}[htb] \epsfysize = 2.2 in
\centerline{\epsffile{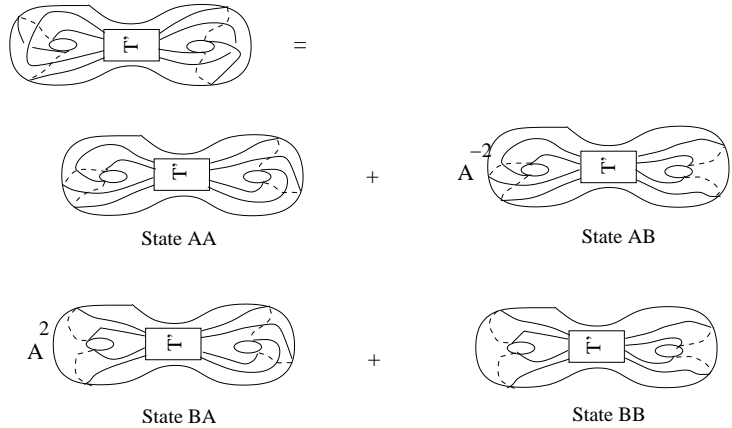}}
\caption{Expansion of a Representation of $ K'_v $ }
\label{fig:fourstate}
\end{figure}

We will refer to the labels assigned to each state later in this section.
We consider a representation of $ K'_v $ in the connected sum of 
two tori. Applying the skein relation to the isolated crossings, 
we obtain an equation with four states. This equation is
illustrated in Figure \ref{fig:fourstate}. 

For some virtual knot diagrams it is possible to determine (or bound)
the virtual genus of the representation without a full expansion of $ T' $.

\begin{figure}[htb] \epsfysize = 0.7 in
\centerline{\epsffile{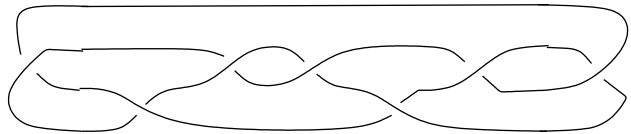}}
\caption{Knot $ K $ with Unknotting Number 2}
\label{fig:fullconcknot}
\end{figure}

We introduce an example constructed from a classical knot diagram with
unknotting number two.  Consider the classical knot diagram $ K' $
with unknotting number 2, shown in Figure \ref{fig:fullconcknot}. The
diagram $K$ has an associated virtualized diagram $ K'_v $ constructed
as above.
\begin{thm}The virtual knot diagram $ K'_v $
has virtual genus two. 
\end{thm}

\begin{figure}[htb] \epsfysize = 0.7 in
\centerline{\epsffile{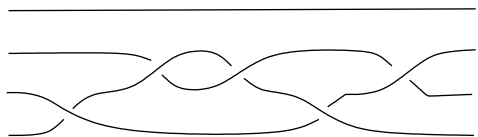}}
\caption{The tangle $T' $}
\label{fig:concknot}
\end{figure}

\textbf{Proof}\qua We decompose $ K $ into two isolated crossings and a classical 
4-4 tangle $ T' $, illustrated in Figure \ref{fig:concknot}.

\begin{figure}[htb] \epsfxsize\hsize
\centerline{\epsffile{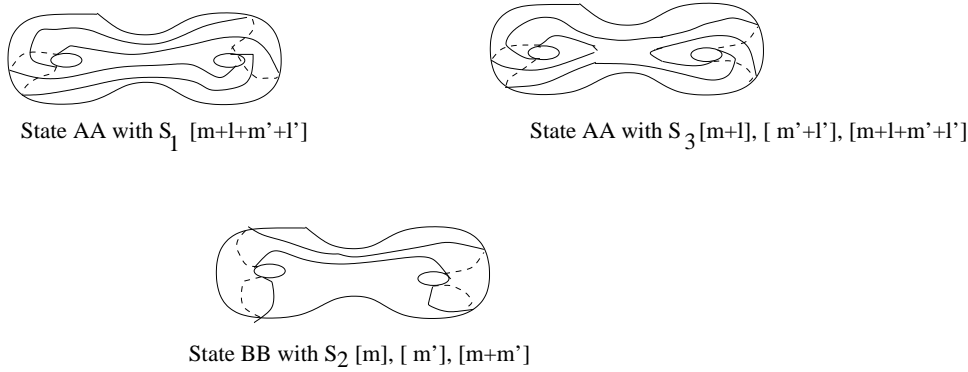}}
\caption{Non-Zero States}
\label{fig:finalstates}
\end{figure}
In Figure \ref{fig:tl4} we list the states obtained from a bracket expansion of a $ 4-4 $ tangle. These states correspond to the generators of the $ 4^{th} $ Temperly-Lieb algebra. 
The states and coefficients obtained from the bracket expansion of
 $ T' $ are:
\begin{gather*} \label{exastates}
A^{-1} s_1 + (A^9 - 2A^5+ 2A ) s_3 + (-A+2A^{-3} - A^{-7}) s_4  \\ 
+ (A^7 - 2A^3 +2 A^{-1} - A^{-5})s_5 + (-A^3 +A^{-1})s_6 
\end{gather*}  
Inserting this expansion into the relation obtained from the skein relation, as shown in Figure \ref{fig:fourstate},
we have the following non-zero states shown in Figure
\ref{fig:finalstates}.
 
These states are sufficient to prevent the presence of 
any cancellation curves. \endproof

\noindent {\bf Acknowledgments}\qua The views expressed herein are
those of the authors and do not purport to reflect the position of the
United States Military Academy, the Department of the Army, or the
Department of Defense.  Much of this effort was sponsored for the
second author by the Defense Advanced Research Projects Agency (DARPA)
and Air Force Research Laboratory, Air Force Materiel Command, USAF,
under agreement F30602-01-2-05022. The U.S. Government is authorized
to reproduce and distribute reprints for Government purposes
notwithstanding any copyright annotations thereon. The views and
conclusions contained herein are those of the authors and should not
be interpreted as necessarily representing the official policies or
endorsements, either expressed or implied, of the Defense Advanced
Research Projects Agency, the Air Force Research Laboratory, or the
U.S. Government.  It gives the second author great pleasure to
acknowledge support from NSF Grant DMS-0245588. (Copyright 2005.)

\Addresses\recd

\end{document}